\documentclass[11pt, a4paper]{article}

\usepackage[english]{babel}
\usepackage[utf8x]{inputenc}
\usepackage[T1]{fontenc}
\usepackage{xparse,amsmath}
\usepackage{tikz}
\usepackage{amsfonts}
\usepackage{amsthm}
\usepackage{amssymb}
\usepackage{graphicx}
\usepackage[shortlabels]{enumitem}

\usepackage{array}
\usepackage{amstext}
\usepackage[shortlabels]{enumitem}
\newcolumntype{C}{>{$}l<{$}}
\usepackage{booktabs}

\usepackage[a4paper,top=3cm,bottom=2cm,left=3cm,right=3cm,marginparwidth=1.75cm]{geometry}

\usepackage{amsfonts}
\usepackage[colorinlistoftodos]{todonotes}
\usepackage{hyperref}
\usepackage{float}
\usepackage{relsize}

\hypersetup{
  colorlinks   = true,    
  urlcolor     = blue,    
  linkcolor    = cyan,    
  citecolor    = red      
}

\title{Domineering games with minimal number of moves}
\author{Rohan Karthikeyan 
\qquad Siddharth Sinha\\
\small Department of Applied Mathematics\\
\small Delhi Technological University\\ 
\small Delhi, India\\
\small\tt fromrohank07@gmail.com\\}

\begin{document}
\theoremstyle{plain}
\newtheorem{theorem}{Theorem}[section]
\newtheorem{problem}[theorem]{Problem}
\newtheorem{definition}[theorem]{Definition}
\newtheorem{remark}[theorem]{Remark}
\newtheorem{observation}[theorem]{Observation}

\maketitle

\begin{abstract}
Domineering is a two-player game played on a checkerboard in which one player places dominoes vertically, while the other places them horizontally. In this paper, we find out the minimum number of moves for a game of domineering to end on several rectangular $m \times n$ boards. We also formulate two problems pertaining to patterns found in the obtained results.

\end{abstract}

\section{Introduction}
Combinatorial games are two-player games with perfect information and no chance moves such as rolling dice or shuffling cards. Some other examples include Chess and Go. The rules are defined in such a way that play will always come to an end with the winner being decided by the game's final move. In \textit{normal play convention}, the first player unable to move loses. In the \textit{misère play convention}, the last player to move loses. \newline

In this paper, we consider the game \textbf{Domineering}, invented by G\"{o}ran Andersson around 1973 and popularized by Martin Gardner \cite{gard}. The version introduced by Andersson and Gardner was the $8 \times 8$ board. Play consists of the two players alternately placing a $1 \times 2$ tile (domino) on adjacent empty squares; Left places vertically and Right places horizontally. The game ends when the player whose turn it is cannot place a piece; the player who cannot place loses - this is the \textit{normal play} condition. Since the board is gradually filled, Domineering is a converging game, the game always ends, and ties are impossible. For more information we would like to direct the attention of the interested reader to the books of Berlekamp, Conway and Guy \cite{ww} and Albert, Nowakowski and Wolfe \cite{albert}. \newline

After stating Huntemann and McKay's \cite{cdp} results of counting Domineering positions satisfying certain properties, we obtain values for the minimum number of moves for a $m\times n$ Domineering game to end in Section \hyperref[sec:Hosanna]{3}, before finishing with further research directions.

\section{Counting Domineering positions}
This section is adapted from \cite{cdp}. We are interested in enumerating positions at the end of the game. Therefore, one defines the following terms:
\begin{definition}
A \textbf{Right end} is a position in which Left potentially has moves available but Right has no moves; a \textbf{Left end} is defined similarly.
\end{definition}
\begin{definition}
A \textbf{maximal position} is a position which is both a left end and a right end - a position in which no player can place a domino.
\end{definition}

\subsection{Maximal Positions}
For counting the maximal domineering positions on a $m \times n$ board, we find the generating function $$F_{m, n}(x, y) = \sum f(a, b)x^a y^b$$ where $f(a, b)$ is the number of maximal domineering positions with $a$ vertical (placed by Left) and $b$ horizontal (placed by Right) dominoes.

\begin{theorem} [Huntemann and McKay \cite{cdp}]
The generating function for the maximal position of an $m \times n$ Domineering board is $$F_{m, n}(x, y) = \mathlarger{\mathlarger{\sum}}\limits_{u \, \in \{0, 1\}^n} \bigg (M_{0, n} + M'_{0, n} \bigg)^m \, (1 + \sum\limits_{i=1}^{n}u_i3^i, 1)$$ where $M_{0, 0} = [1], M_{1, 0} = [0], M_{2, 0} = [0], M'_{0, 0} = [0], M'_{1, 0} = [1], M'_{2, 0} = [0],$ 
\newline
$$M_{0,\, (q+1)} = \begin{bmatrix} M_{2, q} & M_{2, q} & xM_{0, q}\\ M_{1, q} & \mathbf{0} & \mathbf{0}\\ M_{0, q} & M_{0, q} & \mathbf{0} \end{bmatrix},\, \, 
M'_{0,\, (q+1)} = \begin{bmatrix} M'_{2, q} & M'_{2, q} & xM'_{0, q}\\ M'_{1, q} & \mathbf{0} & \mathbf{0}\\ M'_{0, q} & M'_{0, q} & \mathbf{0} \end{bmatrix}$$

$$M_{1,\, (q+1)} = \begin{bmatrix} M_{2, q} & M_{2, q} & xM_{0, q}\\ \mathbf{0} & \mathbf{0} & \mathbf{0}\\ M_{0, q} & M_{0, q} & \mathbf{0} \end{bmatrix},\, \,
M'_{1,\, (q+1)} = \begin{bmatrix} M'_{2, q} & M'_{2, q} & xM'_{0, q}\\ \mathbf{0} & \mathbf{0} & \mathbf{0}\\ M'_{0, q} & M'_{0, q} & \mathbf{0} \end{bmatrix}$$

$$M_{2,\, (q+1)} = \begin{bmatrix} yM_{0, q} & yM_{0, q} & \mathbf{0}\\ \mathbf{0} & \mathbf{0} & \mathbf{0}\\ \mathbf{0} & \mathbf{0} & \mathbf{0} \end{bmatrix}, \, \text{and}\, \, 
M'_{2,\, (q+1)} = \begin{bmatrix} yM'_{0, q} & yM'_{0, q} & \mathbf{0}\\ \mathbf{0} & \mathbf{0} & \mathbf{0}\\ \mathbf{0} & \mathbf{0} & \mathbf{0} \end{bmatrix}$$
\end{theorem}

\subsection{Left and Right Ends}
For counting the Domineering Right ends on a $m \times n$ board, we find the generating function $$R_{m, n}(x, y) = \sum r(a, b)x^a y^b$$ where $r(a, b)$ is the number of Right end positions with $a$ vertical and $b$ horizontal dominoes. The generating polynomial for Left end positions on an $m\times n$ board, which we denote by $L_{m, n}(x, y)$, can be found by obtaining $R_{n, m}(x, y)$ and then switching $x$ and $y$.

\begin{theorem} [Huntemann and McKay \cite{cdp}]
The generating polynomial of Domineering Right ends on an $m\times n$ board is the $(1, 1)$ entry of $(R_{0, n} + R'_{0, n})^m$, denoted by $R_{m, n}(x, y)$, where $R_{0, 0} = [1], R_{1, 0} = [0], R_{2, 0} = [0], R'_{0, 0} = [0], R'_{1, 0} = [1], R'_{2, 0} = [0],$ 
\newline
$$R_{0,\, (q+1)} = \begin{bmatrix} R_{1, q} + R_{2, q} & xR_{0, q}\\ R_{0, q} & \mathbf{0} \end{bmatrix},\, \, 
R'_{0,\, (q+1)} = \begin{bmatrix} R'_{1, q} + R'_{2, q} & xR'_{0, q}\\ R'_{0, q} & \mathbf{0} \end{bmatrix}$$

$$R_{1,\, (q+1)} = \begin{bmatrix} R_{2, q} & xR_{0, q}\\ R_{0, q} & \mathbf{0} \end{bmatrix},\, \, 
R'_{1,\, (q+1)} = \begin{bmatrix} R'_{2, q} & xR'_{0, q}\\ R'_{0, q} & \mathbf{0} \end{bmatrix}$$

$$R_{2,\, (q+1)} = \begin{bmatrix} yR_{0, q} & \mathbf{0}\\ \mathbf{0} & \mathbf{0} \end{bmatrix}, \text{and}\, \, 
R'_{2,\, (q+1)} = \begin{bmatrix} yR'_{0, q} & \mathbf{0}\\ \mathbf{0} & \mathbf{0} \end{bmatrix}$$
\end{theorem}

\section{Results}
\label{sec:Hosanna}
In this section, we state and discuss the results we have obtained and pose some problems not considered in the literature thus far.\newline

One observes that the presence of a monomial of the form $x^a y^b$ in either the Left end ($L_{m, n}(x, y)$) or Right end ($R_{m, n}(x, y)$) or the Maximal end ($F_{m, n}(x, y)$) polynomials indicates the tiling of the $m\times n$ Domineering board with $a$ vertical (placed by Left) tiles and $b$ horizontal (placed by Right) tiles satisfying the required properties of the position.\newline

Now, one might be tempted to think that the minimum number of moves for a $m\times n$ Domineering game to end is the lowest sum of the degrees of $x$ and $y$ in the Left end, Right end or Maximal end polynomials. But, in alternate play, we need to assert that the number of vertical dominoes $a$ and the number of horizontal dominoes $b$ satisfies: $$|a - b| \leq 1$$

Additionally, to get Left ends from $L_{m,n}$ when forcing alternating play, it is not enough to restrict to those terms with $|a-b|\leq 1$. For a position to be a Left end in alternating play, Right must have made the last move. So, we would need $a - b = 0$ or $a - b = -1$, depending on which player started the game. Similarly, for Right ends, we need to enforce the condition that $a - b = 1$ or $a - b = 0$.\newline

Furthermore, we do not need to consider the maximal positions. Recall that all maximal positions are Left ends and Right ends. So, on the one hand we have that the maximal positions are unlikely to be representing shortest play. But on the other hand, even if they are, those positions then appear in the Left end and Right end polynomials anyways.\newline

With these conditions in place, we find the minimum number of moves for a $m \times n$ domineering game to end, which we denote by $\alpha_{m, n}$ as:
\begin{multline}
     \alpha_{m, n} = \min\{a + b \, : x^a y^b \text{occurs in} \, L_{m, n}(x, y) \text{ with} \, a - b \in \{0, -1\} \, \\ \text{(or)}\, R_{m, n}(x, y) \text{ with} \, a - b \in \{0, 1\}\}
\end{multline}

We have incorporated the matrix recurrence relations for the Left and Right ends defined previously into our program written in Maple. The program works by first determining the terms in $L_{m, n}(x, y)$ and $R_{m, n}(x, y)$ where the difference between the powers of $x$ and $y$ is at most $1$. That is because these polynomials actually count all legal Domineering positions satisfying the required properties. Then, we write another procedure to determine the polynomial(s) corresponding to $\alpha_{m, n}$. Additionally, the polynomials were analysed to check whether the term(s) occur in the Left and/or Right end positions.\newline

We have generated results for all rectangular boards up to size $8\times 8$ and some other rectangular boards. Beyond that point, the computations quickly grows too large to fit in the average home computer's main memory. An overview of the results is given in Table \ref{tab:result}. Some preliminary observations are in order:

\begin{table}[ht]
\centering
    \begin{tabular}{C C C C C C C C C C C}
    \toprule
    & \multicolumn{10}{c}{n} \\
    \cmidrule(lr){2-11}
    m & 1 & 2 & 3 & 4 & 5 & 6 & 7 & 8 & 9 & 10 \\
    \midrule
    1 & 0 & 0 & 0 & 0 & 0 & 0 & 0 & 0 & 0 & 0 \\
    2 & 0 & 1^{\textsc{LR}} & 1^{\textsc{R}} & 3^{\textsc{L}} & 3^{\textsc{L}} & 4^{\textsc{LR}} & 4^{\textsc{R}} & 4^{\textsc{R}} & 5^{\textsc{R}} & 5^{\textsc{R}} \\
    3 & 0 & 1^{\textsc{L}} & 2^{\textsc{LR}} & 3^{\textsc{R}} & 3^{\textsc{L}} & 4^{\textsc{L}} & 5^{\textsc{R}} & 5^{\textsc{L}} & 6^{\textsc{L}} & 8^{\textsc{L}} \\
    4 & 0 & 3^{\textsc{R}} & 3^{\textsc{L}} & 4^{\textsc{LR}} & 5^{\textsc{R}} & 7^{\textsc{LR}} & 7^{\textsc{L}} & 8^{\textsc{LR}} & 9^{\textsc{R}} &  \\
    5 & 0 & 3^{\textsc{R}} & 3^{\textsc{R}} & 5^{\textsc{L}} & 7^{\textsc{LR}} & 8^{\textsc{LR}} & 9^{\textsc{LR}} & 11^{\textsc{LR}} &  &  \\
    6 & 0 & 4^{\textsc{LR}} & 4^{\textsc{R}} & 7^{\textsc{LR}} & 8^{\textsc{LR}} & 9^{\textsc{LR}} & 11^{\textsc{LR}} & 12^{\textsc{LR}} &  &  \\
    7 & 0 & 4^{\textsc{L}} & 5^{\textsc{L}} & 7^{\textsc{R}} & 9^{\textsc{LR}} & 11^{\textsc{LR}} & 12^{\textsc{LR}} & 15^{\textsc{LR}} &  &  \\
    8 & 0 & 4^{\textsc{L}} & 5^{\textsc{R}} & 8^{\textsc{LR}} & 11^{\textsc{LR}} & 12^{\textsc{LR}} & 15^{\textsc{LR}} & 16^{\textsc{LR}} &  &  \\
    9 & 0 & 5^{\textsc{L}} & 6^{\textsc{R}} & 9^{\textsc{L}} &  &  &  &  &  &  \\
    10 & 0 & 5^{\textsc{L}} & 8^{\textsc{R}} &  &  &  &  &  &  &  \\
    \bottomrule
    \end{tabular}
    \caption{The minimum number of moves for $m\times n$ Domineering games to end. The superscript over the numbers denote the polynomials (Left end and/or Right end) in which these numbers occur. For instance, the $2\times 6$ game ends in $3$ moves and this position occurs in both the Left end and Right end polynomials. Hence, the entry is denoted $3^{\textsc{LR}}$.}
    \label{tab:result}
\end{table}

\begin{observation}
Note that games on all $m\times 1$ and $1\times n$ boards with $m, n \geq 1$ need 0 minimum moves to end. For example, if we consider the $1 \times n$ boards with Left moving first, Left cannot make a move, and thus the game ends without any moves. The case for the $m \times 1$ boards is symmetrical.
\end{observation}
\begin{observation}
The minimum number of moves for a $m\times n$ Domineering game to end is the same as that for a $n\times m$ Domineering game to end, that is, $\alpha_{m, n} = \alpha_{n, m}$. This is due to the fact that the $n \times m$ board is the negative of the $m\times n$ board. So, this changes only the roles of the two players and not the number of moves.
\end{observation}

On close observation, one can notice the following patterns which, to the best of our knowledge, have not been widely studied:
\begin{observation}
It seems that the minimum number of moves for a $2\times n$ (also, symmetrically the $m \times 2$ boards) Domineering game to end is given by the OEIS \cite{oeis} sequence \href{http://oeis.org/A319198}{A319198}.
\end{observation}
This is also bolstered by the result: $\alpha_{2, 11} = 7^{\textsc{L}}$ (not included in Table \ref{tab:result}). We construct Domineering positions that use the minimum number of moves shown in Table \ref{tab:result} on some $2\times n$ boards.

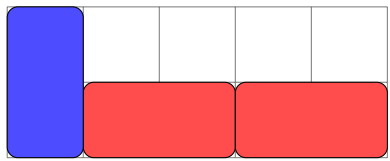
\begin{figure}[H]
\centering
\begin{tikzpicture}
\draw[step=1cm,gray] (-2, -2) grid (3, 0);
\filldraw[fill=red!70, rounded corners] (-1, -2) rectangle (1, -1) {};
\filldraw[fill=blue!70, rounded corners] (-2, -2) rectangle (-1, 0) {};
\filldraw[fill=red!70, rounded corners] (1, -2) rectangle (3, -1) {};
\end{tikzpicture}
\caption{$2 \times 5$ game}
\end{figure}

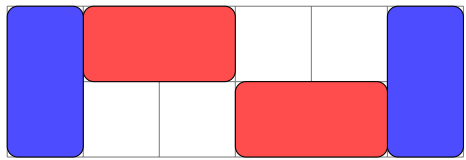
\begin{figure}[H]
\centering
\begin{tikzpicture}
\draw[step=1cm,gray] (-2, -2) grid (4, 0);
\filldraw[fill=red!70, rounded corners] (-1, -1) rectangle (1, 0) {};
\filldraw[fill=blue!70, rounded corners] (-2, -2) rectangle (-1, 0) {};
\filldraw[fill=blue!70, rounded corners] (3, -2) rectangle (4, 0) {};
\filldraw[fill=red!70, rounded corners] (1, -2) rectangle (3, -1) {};
\end{tikzpicture}
\caption{$2 \times 6$ game}
\end{figure}

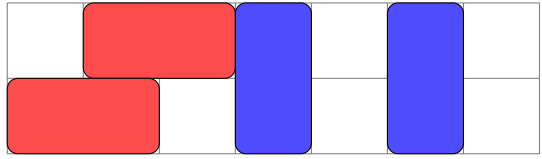
\begin{figure}[H]
\centering
\begin{tikzpicture}
\draw[step=1cm,gray] (-2, -2) grid (5, 0);
\filldraw[fill=red!70, rounded corners] (-1, -1) rectangle (1, 0) {};
\filldraw[fill=blue!70, rounded corners] (1, -2) rectangle (2, 0) {};
\filldraw[fill=blue!70, rounded corners] (3, -2) rectangle (4, 0) {};
\filldraw[fill=red!70, rounded corners] (-2, -2) rectangle (0, -1) {};
\end{tikzpicture}
\caption{$2 \times 7$ game}
\end{figure}

\begin{figure}[H]
\centering
\begin{tikzpicture}
\draw[step=1cm,gray] (-2, -2) grid (6, 0);
\filldraw[fill=red!70, rounded corners] (-1, -2) rectangle (1, -1) {};
\filldraw[fill=blue!70, rounded corners] (2, -2) rectangle (3, 0) {};
\filldraw[fill=blue!70, rounded corners] (4, -2) rectangle (5, 0) {};
\filldraw[fill=red!70, rounded corners] (-1, -1) rectangle (1, 0) {};
\end{tikzpicture}
\caption{$2 \times 8$ game}
\end{figure}
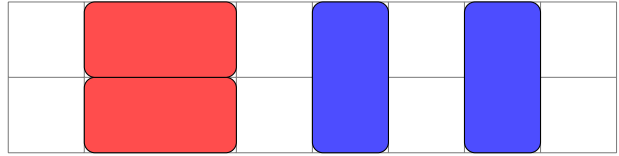

Other sequences of $\alpha_{m, n}$ do not currently appear in the OEIS. Nevertheless, we construct Domineering positions that use the least number of moves on some $3\times n$ boards.

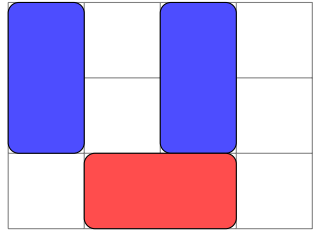
\begin{figure}[H]
\centering
\begin{tikzpicture}
\draw[step=1cm,gray] (-2, -2) grid (2, 1);
\filldraw[fill=blue!70, rounded corners] (-2, -1) rectangle (-1, 1) {};
\filldraw[fill=blue!70, rounded corners] (0, -1) rectangle (1, 1) {};
\filldraw[fill=red!70, rounded corners] (-1, -2) rectangle (1, -1) {};
\end{tikzpicture}
\caption{$3 \times 4$ game}
\end{figure}

\begin{figure}[H]
\centering
\begin{tikzpicture}
\draw[step=1cm,gray] (-2, -2) grid (3, 1);
\filldraw[fill=red!70, rounded corners] (1, -1) rectangle (3, 0) {};
\filldraw[fill=blue!70, rounded corners] (0, -2) rectangle (1, 0) {};
\filldraw[fill=red!70, rounded corners] (-2, -1) rectangle (0, 0) {};
\end{tikzpicture}
\caption{$3 \times 5$ game}
\end{figure}
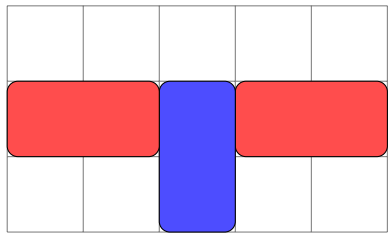

\begin{figure}[H]
\centering
\begin{tikzpicture}
\draw[step=1cm,gray] (-2, -2) grid (4, 1);
\filldraw[fill=red!70, rounded corners] (0, -1) rectangle (2, 0) {};
\filldraw[fill=blue!70, rounded corners] (2, -1) rectangle (3, 1) {};
\filldraw[fill=blue!70, rounded corners] (3, -2) rectangle (4, 0) {};
\filldraw[fill=red!70, rounded corners] (-2, -1) rectangle (0, 0) {};
\end{tikzpicture}
\caption{$3 \times 6$ game}
\end{figure}
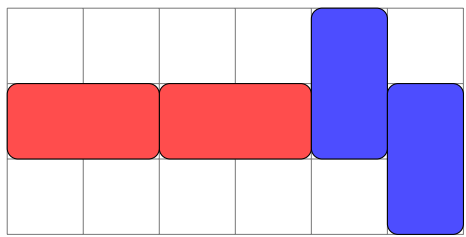

\begin{problem}
$\alpha_{m, n}$ coincides with the minimum number of moves obtained in perfect play only for small $m, n$.
\end{problem}
Particularly, this result holds true only for the $2 \times 2$, $2\times 3$ (and $3\times 2$) and $2 \times 4$ (and $4\times 2$) boards. It's easy to see why this holds true for the $2\times 2$ game. We show the positions of the 2 other boards:

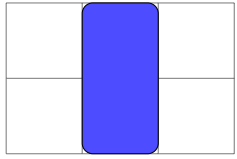
\begin{figure}[H]
\centering
\begin{tikzpicture}
\draw[step=1cm,gray] (-2, -2) grid (1, 0);
\filldraw[fill=blue!70, rounded corners] (-1, -2) rectangle (0, 0) {};
\end{tikzpicture}
\caption{$2 \times 3$ game}
\end{figure}

\begin{figure}[H]
\centering
\begin{tikzpicture}
\draw[step=1cm,gray] (-2, -2) grid (2, 0);
\filldraw[fill=red!70, rounded corners] (-1, -2) rectangle (1, -1) {};
\filldraw[fill=blue!70, rounded corners] (1, -2) rectangle (2, 0) {};
\filldraw[fill=red!70, rounded corners] (-2, -1) rectangle (0, 0) {};
\end{tikzpicture}
\caption{$2 \times 4$ game}
\end{figure}
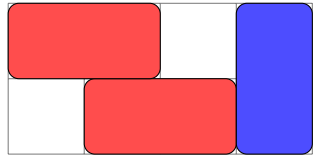

\begin{problem}
For $m \geq 5, n \geq 5$, it seems that $\alpha_{m, n}$ occurs in both the Left end and the Right end polynomials.
\end{problem}

\section{Conclusions and Future research}
We have analysed the minimum number of moves for a $m\times n$ domineering game to end and formulated two new problems. We find two avenues to expand upon the research presented in this paper. One is to continue finding values for unsolved domineering boards using better programs and/or more computer memory. Secondly, one can find the minimum number of moves for a $m\times n$ domineering game to end if both players play perfectly.

\section{Acknowledgements}
We would like to acknowledge Dr. Chandra Prakash Singh, our thesis advisor, for his constant encouragement in the preparation of this manuscript. We would also like to thank Dr. Samrith Ram of IIIT Delhi for introducing us to this wonderful topic of combinatorial game theory. Last but not the least, we would like to thank our parents for providing moral support all through our lives.


\addcontentsline{toc}{section}{References}
\begin{thebibliography}{1}
\bibitem{albert}
M.K. Albert, R.J. Nowakowski, D. Wolfe.
\textit{Lessons in Play: An Introduction to Combinatorial Game Theory},
A K Peters, 2007.

\bibitem{ww}
E. Berlekamp, J.H. Conway, R.K. Guy.
\textit{Winning Ways for your Mathematical Plays},
vol. 1 (2001), vols. 2, 3 (2003), vol. 4 (2004). A K Peters.

\bibitem{gard}
M. Gardner.
\textit{Mathematical games: cram, crosscram and quadraphage: new games having elusive winning strategies}.
Sci. Am., \textbf{230} (1974), 106-108.

\bibitem{cdp}
S. Huntemann, N. A. McKay.
\textit{Counting Domineering positions}.
arXiv preprint, {\tt arXiv:1909.12419} (2019). 

\bibitem{oeis}
N. J. A. Sloane et. al., The On-Line Encyclopedia of Integer Sequences,
https://oeis.org, 2020.\newline
\end{thebibliography}
\end{document}